\documentclass{amsart}
\usepackage{mathrsfs}
\usepackage{amssymb}

\newtheorem{theorem}{Theorem}[section]

\newtheorem{corollary}[theorem]{Corollary}
\newtheorem{proposition}[theorem]{Proposition}
\newtheorem{definition}[theorem]{Definition}
\newtheorem{remark}[theorem]{Remark}

\newtheorem{hypotheses}[theorem]{Hypotheses}

\DeclareMathOperator{\SW}{\mathcal {SW}}
\DeclareMathOperator{\Spec}{Sp}
\DeclareMathOperator{\Hom}{Hom}
\DeclareMathOperator{\RHom}{RHom}
\DeclareMathOperator{\Rhom}{RHom}

\DeclareMathOperator{\map}{map}
\DeclareMathOperator{\Ho}{Ho}
\DeclareMathOperator{\id}{id}
\DeclareMathOperator{\op}{op}
\DeclareMathOperator{\HoC}{\Ho(\C)}

\DeclareMathOperator{\ev}{Ev}

\newcommand{\mZ}{{\mathbb Z}}

\newcommand{\mS}{{\mathbb S}}
\newcommand{\mI}{{\mathbb I}}
\newcommand{\mIo}{{\mathbb I}^{-1}_{c}}
\newcommand{\mIn}{{\mathbb I}^{-n}_{c}}
\newcommand{\mInm}{{\mathbb I}^{-(n+m)}_{c}}
\newcommand{\mIm}{{\mathbb I}^{-m}_{c}}
\newcommand{\onI}{\omega^n{\mathbb I}}
\newcommand{\omI}{\omega^m{\mathbb I}}
\newcommand{\ompnI}{\omega^{m+n}{\mathbb I}}
\newcommand{\ooI}{\omega^1{\mathbb I}}
\newcommand{\ozI}{\omega^0{\mathbb I}}
\newcommand{\C}{{\mathcal C}}
\newcommand{\D}{{\mathcal D}}
\newcommand{\T}{{\mathcal T}}
\newcommand{\W}{{\mathscr W}}
\newcommand{\spec}{Sp^{\Sigma}}

\newcommand{\ssets}{\mathcal{S}_*}
\newcommand{\iso}{\cong}
\newcommand{\bd}{\bf{\cdot}}
\newcommand{\rscript}[1]{^{\mbox{\scriptsize {#1}}}}

\newcommand{\sm}{\wedge}
\newcommand{\Smash}{\wedge}

\newcommand{\boxprod}{\mathbin{\square }}

\newcommand{\mc}{\colon \,}

\renewcommand{\to}{\longrightarrow}
\newcommand{\varrow}[1]{\hbox to #1{\rightarrowfill}}
\newcommand{\varl}[2]{\stackrel{#2}{\hbox to #1{\leftarrowfill}}}
\newcommand{\varrx}[1]{\stackrel{#1}{\hbox to {1cm}{\rightarrowfill}}}
\newcommand{\varr}[1]{\stackrel{#1}{\longrightarrow}}
\newcommand{\parallelarrows}{\begin{array}{c} {
{\longrightarrow}}  \vspace{-0.35cm} \\ {
{\longrightarrow}} \end{array}}

\begin{document}

\title{Monoidal uniqueness of stable homotopy theory}
\date{\today; 2000 AMS Math.\ Subj.\ Class.: 55U35, 55P42}
\author{Brooke Shipley}
\thanks{Research partially supported by an NSF grant.}
\address{Department of Mathematics \\ University of Chicago\\  
Chicago, IL 60637 and \\} 
\address{Department of Mathematics \\ Purdue University \\  
West Lafayette, IN 47907} 
\email{bshipley@math.uchicago.edu,  bshipley@math.purdue.edu}

\begin{abstract}
We show that the monoidal product on the stable homotopy category
of spectra is essentially unique.  This strengthens work of this
author with Schwede on the uniqueness of models of the stable homotopy
theory of spectra.  As an application we show that with an added assumption 
about underlying model structures Margolis' axioms uniquely determine 
the stable homotopy category of spectra up to monoidal equivalence.   
Also, the equivalences constructed here give a unified construction of 
the known equivalences of the various symmetric monoidal categories of 
spectra ($S$-modules, $\W$-spaces, orthogonal spectra, simplicial functors) 
with symmetric spectra.  The equivalences of modules, algebras and commutative 
algebras in these categories are also considered. 
\end{abstract}

\maketitle

\section{Introduction}

The homotopy category of spectra,
obtained by inverting the weak equivalences of spectra, 
has long been known to have a symmetric monoidal product (or tensor
product) induced by the smash product~\cite{boardman-spectra, vogt}.  
Recently, several categories of spectra have
been constructed which have symmetric monoidal smash products 
even before the weak equivalences are inverted~\cite{ekmm, hss, lydakis, mmss}.
Such categories are of interest because they facilitate the development of 
algebraic constructions such as ring spectra and module spectra.  
In each of these examples, inverting the
weak equivalences recovers the standard homotopy category of spectra with the 
standard smash product.  This raises the question of whether this is
forced.  In this paper we consider this question about 
the uniqueness properties of the monoidal product on categories 
of spectra and on the homotopy category of spectra.

Each of these categories of spectra is in fact a highly structured
category.  This structure includes a {\em simplicial Quillen model structure} 
which
encodes standard homotopy theoretic constructions~\cite[Chapter II \S 2]{Q}.  
The symmetric monoidal
product is also compatible with this model structure so that the
derived product induces a symmetric monoidal product on the homotopy  category 
where the weak equivalences have been inverted.  
Another common property of these
model categories of spectra is that they are stable - the suspension
functor is invertible up to homotopy (with inverse the loop functor);
see Definition~\ref{def-stable model category}. 
A category with such compatible structures and a cofibrant 
desuspension of the unit is called a 
{\em stable simplicial monoidal model category}; see Definition~\ref{def-smmc}.
(A cofibrant desuspension of the unit exists in these categories of spectra and
any model category for which the conclusion of Theorem~\ref{thm-intro-simp}
below holds; see Remarks~\ref{rem-simp-de} and~\ref{rem-simp-desuspension}.) 

Instead of restricting to these known model categories of spectra we
consider any stable simplicial monoidal model category $\C$.  First
we show that symmetric spectra, $\spec$~\cite{hss}, is initial
among such model categories by constructing a functor from $\spec$ to 
$\C$ which is simplicial, {\em strong monoidal}
and a {\em left Quillen adjoint}.
That is, the functor is compatible with the simplicial action, 
the monoidal product and the model category structure. 
(See~\cite[II \S 2]{Q}, Definition~\ref{def-monoidal-functor} 
and Definition~\ref{def-Quillen pair}).  

\begin{theorem}\label{thm-intro-simp}
There is a simplicial, strong monoidal, left Quillen functor 
from $\spec$ to $\C$ for $\C$ any stable simplicial monoidal model category.   
That is, the positive stable model category on symmetric spectra
is initial among the stable simplicial monoidal model categories. 
\end{theorem}

This statement is proved as Theorem~\ref{thm-simplicial}.  The
positive stable model structure on $\spec$ is recalled
in Definition~\ref{def-psmc}.  This positive model category captures
the same homotopy theory as the standard model category (i.e, they are
{\em Quillen equivalent}, see Definition~\ref{def-Quillen pair}) but is initial
because the sphere spectrum is not positive cofibrant; 
see Theorem~\ref{thm 4.2}.  

Under additional assumptions on $\C$, the functor 
constructed in Theorem~\ref{thm-intro-simp} is a Quillen
equivalence.  Hence these additional assumptions uniquely specify
the models for the homotopy theory of spectra among the stable simplicial 
monoidal model categories. 
The first additional assumption here is that the unit $\mI$
of the monoidal product in $\C$ is a {\em small, weak generator} in the
homotopy category of $\C$.  This is equivalent to asking that 
$[\mI, -]^{\Ho(\C)}$ commutes with coproducts and detects isomorphisms;
see Definition~\ref{def-generator, small}.  For example, the sphere spectrum, 
$\mS$, is a small, weak generator of the homotopy category of spectra.  We also
ask that $[\mI,\mI]^{\Ho(\C)}$ is freely generated as a $\pi^s_*$-module
by the identity map of $\mI$, as holds for $\mS$ in the homotopy category 
of spectra. 

\begin{theorem}\label{cor-intro-simp} Let $\C$ be a stable simplicial
monoidal model category. 
There is a monoidal Quillen equivalence from the positive stable model 
structure on $\spec$ to $\C$ if and only if 
the unit, $\mI$, is a small weak
generator for which $[\mI, \mI]^{\Ho(\C)}_*$ is freely generated as a 
$\pi_*^s$-module by the identity map of $\mI$. 
\end{theorem}

This shows that up to monoidal Quillen equivalence there is a unique 
stable simplicial monoidal model category of spectra which 
satisfies the above hypotheses.
Other equivalent conditions are stated in Theorem~\ref{cor-simp}. 
In Section~\ref{sec-rings} this uniqueness is extended to modules, 
algebras and commutative algebras.
Remark~\ref{rem-simp-desuspension} shows that the monoidal 
Quillen equivalences constructed in Theorems~\ref{thm-simplicial},
~\ref{thm-rings} and ~\ref{thm-ca} recover
and unify those constructed in~\cite{mmss} and~\cite{sch-comparison}
between $S$-modules, orthogonal spectra, $\W$-spaces, simplicial
functors and symmetric spectra and the associated categories of
modules and algebras.  The conditions on the unit in 
Theorem~\ref{cor-intro-simp} were first studied in~\cite{SS3}.
There we considered the uniqueness of model categories of
spectra but ignored the monoidal product structure.  

Theorems~\ref{thm-intro-simp} and~\ref{cor-intro-simp}
give the most highly structured uniqueness properties of the monoidal
product on the model category level.  Next we consider a 
weaker situation which is still strong enough to establish 
uniqueness properties of the
monoidal product on the homotopy category.  This less structured
situation also provides an approach to Margolis' Conjecture,
see Theorem~\ref{thm-margolis} below, with fewer 
hypotheses than would be required using the above statements.  

On the homotopy category, 
Corollary~\ref{cor-hss} shows that under weak hypotheses
if there is a natural transformation $A \sm B \to A\sm' B$ between 
two monoidal products 
which both have the sphere spectrum, $\mS$, as the 
unit, then this transformation is an isomorphism on all objects.  Thus, 
the main obstruction to showing that two monoidal products
are equivalent is constructing a natural transformation between them.

To construct such natural transformations we consider the model categories 
of spectra, rather than the homotopy category.  Here we consider
{\em stable monoidal model categories}, that is, 
stable model categories $\C$ with a compatible monoidal product and
a cofibrant desuspension of the unit; 
see Definition~\ref{def-non-simp}.
We construct a functor from the homotopy category of spectra, $\Ho(\Spec)$,
to the homotopy category of $\C$, $\Ho(\C)$   
which induces a natural isomorphism 
between the smash product on  $\Ho(\Spec)$ and the derived product on 
$\Ho(\C)$ (i.e, a strong monoidal functor,  
see Definition~\ref{def-monoidal-functor}). 

\begin{theorem}\label{thm-intro-non-simp}
Let $\C$ be a stable monoidal model category. 
Then there is a strong monoidal functor from $\Ho(\Spec)$ to $\Ho(\C)$. 
\end{theorem}

This statement is proved as Theorem~\ref{thm-non-simp}.  
As above, with added hypotheses on the unit this strong monoidal functor induces
a structured monoidal equivalence between $\Ho(\C)$ and $\Ho(\Spec)$. 

\begin{theorem}\label{cor-intro-non-simp} 
Let $\C$ be a stable monoidal model category.  There is a 
$\pi_*^s$-linear, triangulated, monoidal equivalence between the
homotopy category of $\C$ and the homotopy category of spectra
if and only if the unit, $\mI$, is a small weak
generator for which $[\mI, \mI]^{\Ho(\C)}_*$ is freely generated as a 
$\pi_*^s$-module by the identity map of $\mI$. 
\end{theorem}

This shows that the only monoidal
product on $\Ho(\Spec)$ which has an underlying model satisfying these
weak hypotheses is the usual smash product.  
This statement and several other equivalent conditions are 
proved as Theorem~\ref{cor-non-simp}.

We apply these results to Margolis' Conjecture from~\cite{margolis}.
Margolis introduced axioms for a {\em stable homotopy category} which
basically ensure that it is the structured completion of the 
Spanier-Whitehead category of finite CW complexes; see 
Definition~\ref{def-margolis}.
He then conjectured that these axioms uniquely determine the stable homotopy
category of spectra.  Here we add the assumption that the category has
an underlying stable monoidal model category; see 
Definition~\ref{def-underlying}. 

\begin{theorem}\label{thm-margolis}
Suppose that $\mathscr{S}$ is a stable homotopy category in the sense 
of~\cite[Chapter\ 2 \S 1]{margolis} which has an underlying stable
monoidal model category.
Then $\mathscr{S}$ is monoidally equivalent to the stable homotopy 
category of spectra.
\end{theorem}

{\em Acknowledgments: } 
The initial impetus for this paper was an observation of Hovey  
which appeared in an early version of~\cite{hss}.  Corollary~\ref{cor-hss} 
is a modification and generalization of that observation. 
This work continues the study begun in~\cite{SS3} on
uniqueness properties of model categories of spectra where the monoidal
product was ignored.  The construction of the functors in the simplicial
case, see Section~\ref{sec-simp}, builds on the special cases developed
in~\cite{mmss} and~\cite{sch-comparison}.  The construction of functors
in the non-simplicial case, see Section~\ref{sec-non-simp}, builds on the 
treatment of cosimplicial
resolutions in~\cite{SS3}.  The new ingredient is a monoidal product on
cosimplicial resolutions that has not been considered before.
I would also like to thank Mike Mandell and Charles
Rezk for helpful suggestions during this project.   

\section{Model category preliminaries} \label{sec-stable}

In this section we recall the relevant definitions.
A monoidal model category is a model category with a compatible
symmetric monoidal product.  Note that we do require the product to be
symmetric even though that term is suppressed in the name `monoidal model
category'.  The compatibility is expressed by the pushout product axiom
below.  This compatibility is analogous to the simplicial axiom of 
~\cite[Chapter II \S 2]{Q}.  In particular, the product on a monoidal model category
induces a derived product on the homotopy category which is symmetric
monoidal.  Monoidal model categories have been studied in~\cite{SS1} and
~\cite{hovey-book}.  Here, instead of requiring a closed monoidal
structure, we use the weaker hypotheses that the product commutes
with colimits.

\begin{definition}\label{def-monoidal} 
{\em A model category $\C$ is a
{\em monoidal model category} if it is endowed
with a symmetric monoidal structure which commutes with colimits and
satisfies the following pushout product axiom and unit axiom.
We denote the symmetric monoidal product by $\sm$ and the unit
by $\mathbb I$.} 

{ Pushout product axiom. }{\em Let $i\mc A \to B$ and
$j\mc X \to Y $ be cofibrations in
$\C$.
Then the map
\[ i \boxprod j\mc A\Smash Y \cup_{A\Smash X} B\Smash X \to  B\Smash Y \]
is a cofibration which is a weak
equivalence if either $i$ or $j$ is a weak equivalence.}

{Unit axiom. }{\em If the unit is not cofibrant then fix a cofibrant
replacement $u\mc Q\mI \to \mI$ which is a trivial fibration from a 
cofibrant object $Q\mI$.  Then for any cofibrant object $X$
the map $u\sm X \mc Q\mI \sm X \to \mI \sm X$ is a weak equivalence.}
\end{definition}

\begin{definition}\label{def-monoidal-functor} 
{\em A functor $F\mc \C \to \D$
between symmetric monoidal categories is {\em lax monoidal} if there
is a map $\eta\mc\mI_{\D} \to F(\mI_{\C})$ and a transformation
$\phi\mc FA \sm_{\D} FB \to F( A \sm_{\C} B)$, natural in both variables,
such that the coherence diagrams for commutativity, associativity and
unital properties commute.  The functor $F$ is {\em strong monoidal}
if $\eta$ and $\phi$ are isomorphisms.   
}\end{definition}

Next we define the appropriate equivalences of model categories and
monoidal model categories.

\begin{definition} \label{def-Quillen pair}
{\em A pair of adjoint functors between model categories is a 
{\em Quillen adjoint pair} if the right adjoint preserves trivial fibrations and
fibrations between fibrant objects.  This is equivalent to the usual
definition~\cite[Definition 1.3.1]{hovey-book} by~\cite[Corollary A.2]{dugger}.
We regard a Quillen adjoint pair as a map of model categories in the
direction of the left adjoint. 
A Quillen adjoint pair induces adjoint total derived functors 
between the homotopy categories~\cite[Chapter I \S 4 Theorem\ 3]{Q}. 
A Quillen functor pair is a {\em Quillen equivalence} if the total 
derived functors are adjoint equivalences of the homotopy categories.
A {\em monoidal Quillen equivalence} is a Quillen equivalence between
monoidal model categories with a strong monoidal left adjoint functor $L$
such that $L(Q\mI) \to L(\mI)$ is a weak equivalence.  An equivalence
of homotopy categories via strong monoidal functors is called a
{\em monoidal equivalence}.    If one functor in an 
adjoint equivalence is strong monoidal then so is the other,  so both the left 
and right total derived functors of a monoidal Quillen equivalence are strong
monoidal.  Hence a monoidal Quillen equivalence induces a monoidal equivalence
on the homotopy categories. 
}\end{definition}

In this paper we actually consider only stable model categories.
Recall from \cite[Chapter I \S 2]{Q} or \cite[Definition 6.1.1]{hovey-book}
that the homotopy category of a pointed
model category supports a suspension functor $\Sigma$ with a
right adjoint loop functor $\Omega$.

\begin{definition} \label{def-stable model category}
{\em  A {\em stable model category} is a pointed, complete and cocomplete 
category with a model category structure for which
the functors $\Omega$ and $\Sigma$ on the homotopy category are
inverse equivalences.}
\end{definition}

Certain extra structures on the homotopy category of a stable model category 
are key here.  The homotopy category is naturally a triangulated 
category~\cite{verdier}. 
The suspension functor defines the shift functor and the cofiber sequences
of~\cite[Chapter I \S 3]{Q} define the distinguished triangles (the fiber sequences
agree up to sign~\cite[Theorem 7.1.11]{hovey-book}); see \cite[Proposition 7.1.6]{hovey-book}
for more details.
We have required a stable model category to have all limits and colimits so 
that its homotopy category has infinite sums and products.
The homotopy category of a stable model
category also has a natural action of the ring $\pi_*^s$ 
of stable homotopy groups of spheres~\cite[Construction 2.3]{SS3}. 
If $F:\C\to\D$ is the left adjoint of a Quillen adjoint pair between
stable model categories,
then the total left derived functor $LF:\HoC\to\Ho(\D)$ of $F$ is 
$\pi_*^s$-linear and an exact functor~\cite[Lemma 6.1]{SS3},~\cite[Proposition 6.4.1]{hovey-book}.  

For objects $A$  and $X$ of  a triangulated category $\T$
we denote by $[A,X]^{\T}_*$ the graded abelian homomorphism group
defined by $[A,X]^{\T}_m=[A[m],X]^{\T}$ for $m\in\mZ$, where 
$A[m]$ is the $m$-fold
shift of $A$. If $\T$ is a $\pi_*^s$-triangulated category, then the
groups $[A,X]^{\T}_*$ form a graded $\pi_*^s$-module.

\begin{definition}\label{def-generator, small}
{\em An object $G$ of a triangulated category $\T$ 
is called a {\em weak generator} 
if it detects isomorphisms; i.e., a map $f:X\to Y$ is an isomorphism 
if and only if it induces an isomorphism between the graded abelian 
homomorphism groups $\T(G,X)_*$ and $\T(G,Y)_*$.  
An object $G$ of $\T$ is {\em small} if for any family of objects 
$\{A_i\}_{i\in I}$ whose coproduct exists the canonical map
\[ \bigoplus_{i\in I} \, \T(G, A_i) \ \varrow{1cm} \ 
\T(G, \coprod_{i\in I} A_i) \]
is an isomorphism.}
\end{definition}

\section{Margolis' uniqueness conjecture} \label{sec-margolis}

In this section we apply our monoidal uniqueness results to Margolis'
conjecture about the uniqueness of the stable homotopy category.
Margolis introduced axioms for a stable homotopy category in~\cite{margolis}.
He then conjectured that these axioms uniquely specify the stable homotopy
category of spectra up to a monoidal, triangulated equivalence of categories.  
In~\cite{SS3}, 
any stable homotopy category satisfying Margolis' axioms and having an
underlying model category was shown to be triangulated equivalent to 
the stable homotopy category of spectra.  Here we strengthen that result
to a monoidal, triangulated equivalence.

First we consider a more general setting than Margolis' stable homotopy
categories.   The following proposition shows that under weak hypotheses
a lax monoidal functor between two monogenic, monoidal, triangulated 
categories is strong monoidal. 

\begin{definition}
{\em A {\em monogenic, monoidal, triangulated category} is a triangulated
category $\T$ with arbitrary coproducts and with a symmetric monoidal, 
bi-exact smash product $\sm$ which commutes with suspensions and coproducts 
such that the unit $\mI$ is a small, weak generator.  
}\end{definition}

\begin{proposition}\label{prop-hss}
Assume $(\T, \sm, \mI)$ and $(\T', \sm', \mI')$ are two monogenic, monoidal,
triangulated categories. 
Suppose that $F \mc \T \to \T'$ is a lax monoidal,  
exact functor with unit map $\eta\mc\mI' \to F(\mI)$ and natural 
transformation $\phi\mc FA \sm' FB \to F( A \sm B)$. 
If $F$ commutes with coproducts, $\eta$ is an isomorphism and
$\phi\mc F\mI \sm' F\mI \to F( \mI \sm \mI)$ is an isomorphism,
then $F$ is strong monoidal. 
\end{proposition}

\begin{proof}
Consider the subcategory of objects $A$ in $\T$ such that $\phi\mc FA \sm' F\mI
\to F(A \sm \mI)$ is an isomorphism.  
By the assumptions on $F$, $\sm$ and $\sm'$, both
source and target commute with triangles and coproducts.
So this subcategory is a localizing subcategory which contains $\mI$.
Since $\mI$ is a small, weak generator it follows that this subcategory
is the whole category.  This follows from~\cite[Theorem 2.3.2]{hps}; see also
\cite[Lemma 2.2.1]{SS2}.  Now fix any $A$ and consider the subcategory of objects 
$B$ in $\T$ such that $\phi \mc FA \sm' FB \to F(A \sm B)$ is an
isomorphism.  Again this is a localizing subcategory which 
contains $\mI$, and hence is the whole category.  Thus, $\phi$ is an
isomorphism for any $A$ and $B$.  
\end{proof}

Since the stable homotopy category of spectra is a monogenic, monoidal,
triangulated category, this gives the following corollary.

\begin{corollary}\label{cor-hss}
Assume that $\sm$ and $\sm'$ are
two monogenic, monoidal, triangulated structures on the homotopy category
of spectra, both with unit the sphere spectrum, $\mS$.  
If the identity functor is lax monoidal and the unit map $\eta$ and
the natural transformation $\phi$ evaluated on the unit are isomorphisms, 
then the identity functor gives a monoidal equivalence between 
these two structures.
\end{corollary}

So the only obstruction to showing that the smash product of spectra
is unique up to monoidal equivalence on the homotopy category is constructing 
a natural transformation between any two monoidal products.  Our solution
is to assume there is an underlying stable monoidal model category 
We state this result for Margolis' stable homotopy categories.

\begin{definition}\label{def-margolis}
{\em 
A {\em stable homotopy category} in the sense of \cite[Chapter \ 2 \S 1]{margolis}
is a monogenic, monoidal, triangulated category $\mathscr{S}$ 
with an exact and strong symmetric monoidal equivalence 
$R:\SW_{\mbox{\scriptsize f}}\to {\mathscr{S}}\rscript{small}$
between the Spanier-Whitehead category of finite CW-complexes 
(\cite{spanier-whitehead}, \cite[Chapter \ 1, \S 2]{margolis}) 
and the full subcategory of small objects in $\mathscr{S}$.
}
\end{definition} 

As shown in~\cite[Section 3]{SS3}, such an equivalence
induces a $\pi_*^s$-linear structure on the triangulated category 
$\mathscr{S}$.  In fact, we could weaken the definition above to only 
require that there is such
an equivalence $R$ with the full subcategory of $\SW$ on the spheres
$S^n$ for $n$ an integer.  

\begin{definition}\label{def-underlying}
{\em A stable homotopy category $\mathscr{S}$ has an {\em underlying stable
monoidal model
category} if there is a monoidal, $\pi_*^s$-linear, exact equivalence
$\Phi\mc \mathscr{S} \to \Ho(\C)$ with $\C$ a stable monoidal model
category; see Definition~\ref{def-non-simp}. 
}\end{definition}

\begin{proof}[Proof of Theorem~\ref{thm-margolis}]
Since $\mathscr{S}$ has an underlying stable monoidal model category, there is an
equivalence $\Phi\mc \mathscr{S} \to \Ho(\C)$ with 
all of the properties mentioned in Definition~\ref{def-underlying}.
Since the properties of a small, weak generator are determined on the
homotopy category level, the image $X\in \HoC$ under $\Phi$ of the unit object 
in $\mathscr{S}$ is a small weak generator of the homotopy category of $\C$.
Because the equivalence $\Phi$ is monoidal and $\pi_*^s$-linear, 
$X$ is isomorphic to the unit and satisfies the hypotheses on the unit in 
Theorem~\ref{cor-intro-non-simp}. 
Thus, the homotopy category of $\C$, and hence also $\mathscr{S}$, 
is monoidally equivalent to the ordinary stable homotopy category 
of spectra.
\end{proof}

\section{Simplicial monoidal model categories}\label{sec-simp}

Here we 
construct a Quillen adjoint pair from  
the positive stable model category on $\spec$ to $\C$.  Then,
under additional hypotheses on the unit, this produces a
monoidal Quillen equivalence from $\spec$ to $\C$.  First we
recall the positive model structure from~\cite[Section 14]{mmss}.  

\begin{definition}\label{def-psmc}
{\em In the positive stable model structure on $\spec$ a map $f$ is a weak
equivalence if it is a stable equivalence,~\cite{hss, mmss}. 
The map $f$ is a positive trivial fibration if $f_n$ is a trivial fibration 
for $n > 0$.
The positive cofibrations and positive fibrations are then determined by the 
respective
right and left lifting properties with respect to the trivial fibrations
and the trivial cofibrations.
}\end{definition}

In~\cite{mmss} only symmetric spectra over topological spaces are considered, 
but the arguments can be easily modified to apply to symmetric spectra
over simplicial sets.  The identity
functor from the usual to the positive structure is a right Quillen functor
since (trivial) fibrations are in particular positive (trivial) fibrations.

\begin{theorem}~\cite[Theorem 14.2, Proposition 14.6]{mmss}\label{thm 4.2}
The positive stable model structure on $\spec$ forms a stable, monoidal model
category.  The identity functor induces a monoidal Quillen equivalence
from the positive stable model structure to the usual stable model
structure on $\spec$.
\end{theorem}

Denote the unit in $\spec$ by $\mS$.  Note $\mS$ is not cofibrant in
the positive stable model category.  To fix its cofibrant replacement for
the unit axiom of the monoidal model category structure,
first consider the $n$th evaluation functor $\ev_n$ on symmetric spectra
which lands in $\Sigma_n$-equivariant spaces. 
For $X$ a $\Sigma_n$-space, the left adjoint $F_n'$ is defined by
$(F_n'X)_k \iso \Sigma_{k} \sm_{\Sigma_{n} 
\times \Sigma_{k-n}} (X \sm S^{k-n})$.  This is a slight variant of the free functor
$F_n$ studied in~\cite{hss}.   Note that $F_1' \iso F_1$ and $F_0'\iso F_0$.
Then define the cofibrant replacement of $\mS$
as the weak equivalence $Q\mS = F_1'S^1 \to F_0S^0 = \mS$
induced by the identity map in level one.    

\begin{proposition}\label{prop-pos-stable}
The fibrant objects in the positive stable model structure are the 
positive $\Omega$-spectra.  That is, $X$ is fibrant if $X_n$ is fibrant
for $n > 0$ and $X_n \to \Omega X_{n+1}$ is a weak equivalence for
$n > 0$.  A map $f$ between positive $\Omega$-spectra is a fibration
if each $f_n$ is a fibration for $n > 0$.
\end{proposition}

\begin{proof}  
The description of the fibrant objects follows from~\cite[Theorem 14.2]{mmss}.
The description of the fibrations follows from the fact that the positive
stable model structure is a localization of the positive level model
structure~\cite[Theorem 14.1]{mmss}.  In a localized model structure the fibrations
between fibrant objects are the fibrations in the original model structure.
So here they are the positive level fibrations.  This statement also follows
from the positive variants of~\cite[Lemma 3.4.12]{hss} or~\cite[Proposition 9.5]{mmss}.
\end{proof}

We now define a stable simplicial monoidal model category.  As mentioned in
the introduction, here we require the following technical hypothesis on
the unit which may not be required in other definitions of stable simplicial
monoidal model categories but is needed here; see 
Remark~\ref{rem-simp-desuspension}.

\begin{definition}\label{def-simp-de}
{\em A {\em cofibrant desuspension of the unit} is a cofibrant object $\mIo$ 
with a weak equivalence $\eta\mc\mIo \otimes S^1 \to \mI$.
}\end{definition}
 
Recall that a monoidal model category is a model category $\C$ with
a symmetric monoidal product that is compatible with the model
structure; see Definition~\ref{def-monoidal}.  
Similarly, a {\em simplicial model category} is a model category
with a compatible action of simplicial sets.  A
simplicial functor is a functor that is compatible with this structure.
See~\cite[Chapter II \S 1, 2]{Q}.  

\begin{definition}\label{def-smmc}
{\em A {\em stable simplicial monoidal model category} is a category $\C$
with a stable, simplicial model structure, a monoidal model structure and
a cofibrant desuspension of the unit such that the  
simplicial action commutes with the monoidal product.  That is, for
$X$, $Y$ in $\C$ and $K$ in $\ssets$ there are natural coherent isomorphisms
$(X \sm Y) \otimes K \iso X \sm (Y \otimes K)$.  
}\end{definition}

\begin{remark}\label{rem-simp-de}  
{\em If the unit $\mI$ in $\C$ is fibrant, then a cofibrant desuspension
exists.  Since $\C$ is stable, there is a cofibrant object $X$ whose suspension
in the homotopy category is isomorphic to $\mI$.  Since $\C$ is simplicial
and $X$ is cofibrant its suspension is modeled by $X \otimes S^1$.
Since $X \otimes S^1$ is cofibrant and $\mI$ is fibrant the isomorphism
in the homotopy category is realized by some weak equivalence in $\C$.

A cofibrant desuspension of the unit exists in every known symmetric
monoidal model category
of spectra.  In the diagram categories of spectra investigated in
~\cite{mmss} and their simplicial analogues~\cite{hss, lydakis}, 
(orthogonal spectra, symmetric spectra, and simplicial
functors or $\W$-spaces) the cofibrant desuspension can be chosen as 
the object denoted $F_1S^0$, with the weak equivalence 
$\eta \mc F_1 S^1 \to F_0 S^0$; 
see~\cite[Definition 1.3, Remark 4.7]{mmss}.  The S-modules of~\cite{ekmm} are all fibrant,
so the previous paragraph applies.}   
\end{remark}

\begin{theorem}\label{thm-simplicial}
Let $\C$ be a stable simplicial monoidal model category. 
Then there exists a Quillen adjoint 
functor pair from the positive stable model structure on $\spec$ to $\C$,  
$\mI \sm - \mc \spec \to \C$ and $\Hom(\mI, -) \mc \C \to \spec$.
These functors are simplicial, the left adjoint 
$\mI \sm -$ is strong monoidal, and $\mI \sm Q\mS \to \mI \sm \mS$ is
a weak equivalence. 
\end{theorem}

Remark~\ref{rem-simp-desuspension} below shows that the existence of such a 
Quillen adjoint pair implies the existence of a cofibrant desuspension.
Adding conditions on the unit in $\C$ shows this Quillen adjoint pair
is a Quillen equivalence.

\begin{theorem}\label{cor-simp}
Let $\C$ be a stable simplicial monoidal model category. 
The following conditions are equivalent:
\begin{enumerate}
\item There is a $\pi^s_*$-linear triangulated equivalence 
from the homotopy category of $\spec$ to the homotopy category of $\C$ 
which takes the unit $\mI$ of the monoidal product in $\C$ to the unit 
$\mS$ of $\spec$.
\item The unit, $\mI$, is a small weak generator for which 
$[\mI, \mI]_*^{\Ho(\C)}$ is freely generated as a $\pi_*^s$-module by the 
identity map of $\mI$.
\item There is a simplicial, monoidal Quillen equivalence from the 
positive stable model structure on $\spec$ to $\C$.   
\item There is a zig-zag of monoidal Quillen equivalences between the 
usual stable model structure on $\spec$ and $\C$. 
\end{enumerate}
\end{theorem}

\begin{proof}
Condition (1) implies condition (2) since the properties of $\mI$ mentioned
in (2) hold for $\mS$ and are determined by the $\pi_*^s$-linear 
triangulated homotopy category.  Condition (3) implies condition (4)
because the positive stable model structure is monoidally Quillen equivalent
to the usual stable model structure on $\spec$ by Theorem~\ref{thm 4.2}.  
Since Quillen functors induce $\pi_*^s$-linear triangulated functors
on the homotopy categories by~\cite[Proposition 6.4.1]{hovey-book} and~\cite[Lemma 6.1]{SS3}
and monoidal functors preserve the unit, condition (4) implies condition (1).

Next we show that given condition (2) the simplicial Quillen adjoint pair 
constructed
in Theorem~\ref{thm-simplicial} is a Quillen equivalence.
Since $\mI \sm -$ is strong monoidal, $\mI \sm \mS \iso \mI$.  
Also, $\mI \sm Q\mS \to \mI \sm \mS$ is a weak equivalence, so
$\mI \sm^L \mS \iso \mI$.  The total derived functor $\mI \sm^L -$ is exact 
by~\cite[Proposition 6.4.1]{hovey-book}.  So $\mI \sm^L \mS[n] \iso \mI[n]$  
where $X[n]$ denotes the $n$th shift of $X$ for any integer $n$. 
This isomorphism and the derived adjunction for $ \mI \sm^L -$ and
$\RHom(\mI, -)$ produce the following natural isomorphisms
\[ \pi_*\RHom(\mI, Y) \iso [\mS[*], \Rhom(\mI, Y)] \iso [\mI, Y]_*^{\HoC}.
\]
Since $\mI$ is a weak generator, $\Rhom(\mI, -)$ detects isomorphisms.
So to show that this pair is a Quillen equivalence we need to show that
for any symmetric spectrum $A$ the 
unit of the adjunction $A \to \Rhom(\mI, \mI \sm^L A)$ is an isomorphism.
Consider the full subcategory $\T$ of such objects. 
For $A = \mS$ in homotopy this map is the map
$[\mS, \mS]_* \to [\mI, \mI]_*$ induced by $\mI \sm^L -$. 
This map of free $\pi^s_*$-modules takes the
identity map of $\mS$ to the identity map of $\mI$.  Hence it is also an 
isomorphism 
by condition (2).  So $\mS$ is contained in $\T$.   Since $\mI$ is small,
$\Rhom(\mI, -)$ commutes with coproducts by the display above.  Hence,
since left adjoints commute with coproducts and total derived functors
between stable model categories are exact,
the composite $\Rhom(\mI, \mI \sm^L -)$ is an exact functor which commutes
with coproducts.  So $\T$ is a localizing subcategory which contains
the generator $\mS$ of symmetric spectra.  Hence $\T$ is the whole category.
Thus, these derived functors induce an equivalence of homotopy categories. 
\end{proof}

\begin{proof}[Proof of Theorem~\ref{thm-simplicial}]
We first construct the functor $\Hom(\mI, -) \mc \C \to \spec$.  
Let $\mIo$ be a cofibrant desuspension of the unit in $\C$
with a weak equivalence $\eta\mc\mIo \otimes S^1 \to \mI$.
Let $\mIn=(\mIo)^{\sm n}$ be the $n$-fold smash product of $\mI$ 
where $X^0= \mI$.  Notice that in general $\mI^0_c$
is not cofibrant.  For $Y$ in $\C$, define the $n$th level of $\Hom(\mI, Y)_n$ 
to be the simplicial mapping space $\map_{\C}(\mIn, Y)$.  The symmetric
group on $n$ letters acts on $\mIn$ by permuting the factors and
hence also acts on $\Hom(\mI, Y)_n$.  The structure map 
\[\map_{\C}(\mIn,Y) \to \Omega^m\map_{\C}(\mInm, Y)\iso 
\map_{\C}(\mInm \otimes S^m, Y)
\]
is induced by $\map_{\C}(\sigma, Y)$ with $\sigma$ defined as 
\[\sigma_{n,m}\mc\mInm \otimes S^m \iso \mIn \sm (\mIo \otimes S^1)^m
\varr{\id \sm (\eta)^m} \mIn \sm (\mI)^m \iso \mIn.\]  Since the
adjoint of $\map_{\C}(\sigma, Y)$ is $\Sigma_n \times \Sigma_m$ equivariant, 
this makes $\Hom(\mI, Y)$ into a symmetric spectrum.  Here we have used the 
fact that the simplicial action and the monoidal product commute. 
$\Hom(\mI, -)$ is an example of a categorical construction described
in~\cite[I.2]{mm}.  

Since $\C$ is a simplicial model category and $\mIn$ is cofibrant for $n > 0$, 
$\Hom(\mI, -)$ of a (trivial) fibration is a (trivial) fibration in
levels $n > 0$.
Since $\mIo \otimes S^1$ is cofibrant $\eta$ factors as
$\mIo \otimes S^1 \to Q\mI \to \mI$ where $Q\mI \to \mI$ is the fixed
cofibrant replacement of $\mI$ given in the monoidal model structure on $\C$.
Since $\C$ is monoidal and $\eta$ is a weak equivalence, $\sigma_{n,1}$
is a weak equivalence between cofibrant objects for $n > 0$. 
Hence $\Hom(\mI, -)$ takes a fibrant object to a positive $\Omega$-spectrum,
which is a fibrant object in the positive stable model structure.
Thus $\Hom(\mI, -)$ takes trivial fibrations to positive
trivial fibrations and fibrations to positive fibrations between
positive fibrant objects by Proposition~\ref{prop-pos-stable}.  
So $\Hom(\mI, -)$ is a right Quillen adjoint.

Next we consider the left adjoint $\mI \sm - \mc \spec \to \C$.
Using the definition of $F_n'X$, $\mI \sm F_n'X$ is isomorphic to  
$\mIn \otimes_{\Sigma_n} X$ since both corepresent the
functor which takes $Y$ in $\C$ to the space of $\Sigma_n$-equivariant
maps from $X$ to $\Hom(\mI, Y)_n$.  So $\mI \sm Q\mS \to \mI \sm \mS$ is
isomorphic to the weak equivalence $\eta \mc \mIo \otimes S^1 \to \mI$.

To evaluate $\mI \sm -$ on an arbitrary symmetric spectrum $A$, note that $A$
can be built as the coequalizer of the following diagram:
\[  \bigvee_n  F_{n+1}'(\Sigma_{n+1} \sm_{\Sigma_n} (A_n \sm S^1))
\parallelarrows \bigvee_n  F_n'A_n\]
Here one map is induced by the map $A_n \sm S^1 \to A_{n+1}$ and
the other is induced by smashing $F_n'A_n$ with the map 
$F_1' S^1 \to F_0'S^0$ which is the adjoint of the identity map on $S^1$
in level one.   Since $\mI \sm -$ must commute with colimits, $\mI \sm A$ is 
defined as the coequalizer of the diagram:
\[  \bigvee_n  \mI_c^{-(n+1)} \otimes_{\Sigma_n} (A_n \sm S^1)
\parallelarrows \bigvee_n  \mIn \otimes_{\Sigma_n} A_n\]
Again the first map is induced by $A_n \sm S^1 \to A_{n+1}$ and the
second map uses the fact that the simplicial action and monoidal product
in $\C$ commute to give the isomorphism 
\[ 
\mI_c^{-(n+1)} \otimes_{\Sigma_n} (A_n \sm S^1) \iso (\mIn \otimes_{\Sigma_n} 
A_n) \sm (\mIo \otimes S^1) \] 
along with the map $\eta \mc \mIo \otimes S^1 \to \mI$.

Next we consider the monoidal properties of these adjoint functors.
First $\Hom(\mI,-)$ is lax monoidal; since the simplicial action and
monoidal product commute the product of maps induces
$\map_{\C}(\mIn, A) \sm \map_{\C}(\mIm, B) \to \map_{\C}( \mInm, A \sm B)$.
These fit together to give a natural map $\Hom(\mI, A) \sm \Hom(\mI, B)
\to \Hom(\mI, A\sm B)$.  The unit map $F_0'S^0=\mS \to \Hom(\mI, \mI)$ is
given by sending the non-base point of $S^0$ to the identity map of $\mI$
in simplicial degree zero of $\Hom(\mI, \mI)_0 = \map_{\C}(\mI, \mI)$.

The left adjoint of a lax monoidal functor is automatically lax comonoidal.
That is, there are structure maps in the opposite direction of a lax
monoidal functor; see the display below.
The adjoint of the unit map is an isomorphism $\mI \sm \mS \to \mI$.
Denote the adjoint pair by $L$ and $R$.  Then the counit and unit
of the adjunction and the lax monoidal structure of $R$ give
\[ L( A \sm B) \to L(RLA \sm RLB) \to LR(LA \sm LB) \to 
LA\sm LB. \]
Here in fact $L= \mI \sm -$ is strong monoidal because this map is an 
isomorphism.  To show this we only need to consider the special case
where $A= F_n'X$ and $B=F_m'Y$ for $X$ a $\Sigma_n$-space and $Y$ a 
$\Sigma_m$-space since the general case follows by using the
coequalizer diagrams above.  Then 
\[  L( A \sm B) = LF_{n+m}' (\Sigma_{n+m} \sm_{\Sigma_n \times \Sigma_m} 
X \sm Y) \iso  \mInm \otimes_{\Sigma_n \times \Sigma_m} X \sm Y.\]  
Again commuting the simplicial action and the monoidal product
shows this last term is isomorphic via the transformation displayed above
to $(\mIn \otimes_{\Sigma_n} X) \sm (\mIm \otimes_{\Sigma_m} Y) = 
LA \sm LB$.  These monoidal properites also follow from the more
general treatment in~\cite[I.2]{mm}.  

Finally, these adjoint functors $\Hom(\mI, -)$ and $\mI \sm -$ are
simplicial functors.  This follows by various adjunctions from the
isomorphism $\Hom(\mI, Y^K) \iso \Hom(\mI, Y)^K$ given by the
simplicial structure on $\C$.    
\end{proof}

\begin{remark}\label{rem-simp-desuspension}
{\em If there is a Quillen adjoint pair from the positive stable
model structure on $\spec$ to $\C$ with a strong monoidal, simplicial left
adjoint $L$ which takes $Q\mS \to \mS$ to a weak equivalence, then
a cofibrant desuspension of the unit exists.
Set $\mIo = L(F_1'S^0)$.
The map $\eta \mc \mIo \otimes S^1 \to \mI$ is then given by 
$L(F_1'S^0) \otimes S^1 \xrightarrow{\varphi} L(F_1'S^1 = Q \mS) \to L(\mS)$
where $\varphi$ is induced by the simplicial structure on $L$.  
The second map is a weak equivalence by assumption.  The first map is a weak
equivalence because it is the cofiber of the weak equivalence $L(F_1S^0)\otimes
\Delta[1]_+ \xrightarrow{\varphi} L(F_1S^0 \otimes \Delta[1]_+)$ by the 
isomorphism $L(F_1S^0) \otimes (S^0 \vee S^0) \xrightarrow{\varphi} 
L(F_1S^0 \otimes (S^0 \vee S^0))$.

This also gives a procedure for recovering the known equivalences between
symmetric monoidal model categories of spectra as 
$\mI \sm -$ and $\Hom(\mI, -)$ for some choice of a cofibrant desuspension
of the unit.  For the monoidal functors constructed in~\cite{mmss} 
($\mathbb {P}$ and $\mathbb{U}$ between orthogonal spectra and symmetric 
spectra and between $\W$-spaces and symmetric spectra), the chosen desuspension
of the unit is $\mathbb{P}(F_1S^0) \iso F_1S^0$~\cite[Definition 1.3, Remark 4.7]{mmss}. 
The monoidal functors ($\Lambda$ and $\Phi$) between $S$-modules and symmetric 
spectra as defined in~\cite{sch-comparison} are 
isomorphic to $\mI \sm -$ and $\Hom(\mI, -)$ with $\Lambda(F_1S^0)\iso
S^{-1}_c$ the chosen desuspension of the unit. 
}
\end{remark}

\begin{remark} {\em If $\C$ is a cofibrantly generated, proper, stable
model category then~\cite[Proposition 4.4]{RSS} shows that $\C$ is 
Quillen equivalent
to a simplicial model category structure on the category of simplicial
objects, $\C^{\Delta^{\op}}$.  If the product on $\C$ commutes with coproducts 
then the level prolongation of the product commutes with the simplicial action.
Using~\cite[Proposition 16.11.1, Theorem 16.4.2]{psh}, one can show that if $\C$ is a monoidal model
category then the simplicial model category from~\cite{RSS} is also monoidal.  
Hence, under these conditions, one can apply the constructions in this
section to the stable simplicial monoidal model category on  
$\C^{\Delta^{\op}}$.  This remark can also be applied if 
$\C$ is simplicial and the product does not commute with the simplicial action
but does commute with coproducts.  We treat the non-simplicial case in even 
more generality in Section~\ref{sec-non-simp}.}
\end{remark}

\section{Modules and Algebras}\label{sec-rings}

In this section, we show that the functors constructed in 
Theorem~\ref{thm-simplicial} induce Quillen adjoint pairs on modules, 
algebras and commutative algebras.  Since
$\mI \sm -$ is strong symmetric monoidal and $\Hom(\mI, -)$ is lax symmetric
monoidal, these
functors restrict to adjoint functors on subcategories of modules and
algebras.   Since we want the restriction of $\Hom(\mI, -)$ to be a right
Quillen adjoint, we assume that in the model structures on categories of
modules or algebras over $\C$ a morphism is a weak equivalence or fibration
if it is one in the underlying model structure on $\C$.  The next proposition
states sufficient conditions for this assumption to hold for modules
and associative algebras.  We treat commutative algebras separately in the 
second part of this section.  

\begin{proposition}\cite[Theorem 4.1]{SS1}
Assume $\C$ is a cofibrantly generated, monoidal model category that
satisfies the monoid axiom~\cite[Definition 3.3]{SS1}.  If the objects in $\C$
satisfy certain smallness conditions~\cite[Lemma 2.3]{SS1}, then the category
of left $R$-modules (for a fixed monoid $R$) and the category of
$R$-algebras (for a fixed commutative monoid $R$) are model categories
with fibrations and weak equivalences determined in $\C$.
\end{proposition}

\begin{theorem}\label{thm-rings}
Let $\C$ be a stable simplicial monoidal model category with a
monoidal Quillen equivalence from $\spec$ to $\C$ (or any equivalent
condition from Theorem~\ref{cor-simp}) such that the conclusions
of the previous proposition hold.   
If
$\mI \sm -$ preserves weak equivalences between stably cofibrant symmetric
spectra, then
$\mI \sm -$ and $\Hom(\mI, -)$ induce a Quillen equivalence 
\begin{enumerate}
\item from the positive stable model category of 
$R$-modules for $R$ a cofibrant symmetric ring spectrum
to $(\mI \sm R)$-modules, and 
\item from the positive stable model category of 
$R$-algebras for $R$ a commutative symmetric
ring spectrum which is cofibrant as a symmetric spectrum 
to $(\mI \sm R)$-algebras. 
\end{enumerate}
These statements also hold with the usual stable model category replacing
the positive one if $\mI$ is cofibrant.
\end{theorem}

Since $\mS$ is cofibrant as a symmetric spectrum the second 
statement implies that the category of symmetric ring spectra, $\mS$-algebras,
and the category of monoids in $\C$, $\mI$-algebras, are Quillen equivalent. 

\begin{remark}\label{rem-hyp-sats}
{\em The hypothesis that $\mI \sm -$ preserves weak equivalences between
stably cofibrant symmetric spectra is satisfied when $\C$ is any one of
the symmetric monoidal model categories of orthogonal spectra, 
$\W$-spaces~\cite{mmss}, simplicial functors~\cite{lydakis}, or 
$S$-modules~\cite{ekmm}.   Since the unit is cofibrant in the first three
cases, this follows from the next proposition.   This holds in the case of
$S$-modules by~\cite[Theorem 3.1]{sch-comparison} and the fact 
that $\Hom(\mI, -)$ detects and preserves weak equivalences.
}
\end{remark}

\begin{proposition}\label{prop-hyp-sats}
If $\mI$ is cofibrant and $\C$ is a stable simplicial monoidal model
category which satisfies any of the equivalent conditions of 
Theorem~\ref{cor-simp}, then $\mI \sm -$ and $\Hom(\mI, -)$ 
form a Quillen equivalence from the usual stable model category
of symmetric spectra to $\C$.  Hence $\mI \sm -$ preserves weak equivalences
between cofibrant symmetric spectra.
\end{proposition}

\begin{proof}
If $\mI$ is cofibrant, then $\Hom(\mI, -)_0= \map (\mI, -)$ also preserves
(trivial) fibrations.  Hence $\Hom(\mI,-)$ is a right Quillen adjoint
functor from $\C$ to the usual stable model category of symmetric spectra.
The statements follow from the
same proof as given in Theorem~\ref{cor-simp}.
\end{proof}

\begin{proof}[Proof of Theorem~\ref{thm-rings}]
Since the (trivial) fibrations in the categories of $(\mI \sm R)$-modules
and $(\mI \sm R)$-algebras are determined on the underlying category,
the restriction of $\Hom(\mI, -)$ in both cases is still a right 
Quillen adjoint functor to
the positive model structure.  Since $\mI$ is assumed to
be a weak generator by condition (2) of Theorem~\ref{cor-simp},
$\Hom(\mI, -)$ preserves and detects weak equivalences.  So 
by~\cite[Lemma 4.1.7]{hss} we only need to show that $\psi_A \mc A \to 
\Hom(\mI, (\mI \sm A)^f)$ is a weak equivalence for $A$ a positive cofibrant 
object in $R$-modules or $R$-algebras where $(\mI \sm A)^f$ is a fibrant 
replacement.  Since fibrations are
determined on the underlying category, a fibrant replacement as a module
or algebra restricts to a fibrant replacement in $\C$.    

Under the given conditions on $R$, 
if $A$ is cofibrant in the positive model category of 
$R$-modules or $R$-algebras then $A$ is cofibrant as a symmetric spectrum.  
By~\cite[Proposition 14.6]{mmss} the identity functor from the
positive stable model structure on $R$-modules to the
usual stable model structure on $R$-modules is a Quillen left adjoint.
So if $A$ is a positive cofibrant
$R$-module then it is a cofibrant $R$-module.  Since $R$ is assumed to be
cofibrant
as a symmetric ring spectrum it is cofibrant as a symmetric spectrum
by~\cite[Theorem 12.1(v)]{mmss}.  Hence, by~\cite[Theorem 12.1(ii)]{mmss}, $A$ is
cofibrant as a symmetric spectrum.  
Again by~\cite[Proposition 14.6]{mmss}, if $A$ is a positive cofibrant $R$-algebra,
then it is a cofibrant $R$-algebra.  Then by~\cite[Theorem 12.1 (ii), (v)]{mmss}
it follows that $A$ is cofibrant as a symmetric spectrum.

We now show that $\psi_B$ is a weak equivalence for $B$ any cofibrant 
symmetric spectrum. 
It then follows that $\mI \sm -$ and $\Hom(\mI, -)$
restrict to Quillen equivalences on the positive stable model categories of 
$R$-modules and
$R$-algebras.  The proof of Theorem~\ref{cor-simp} shows that 
$\psi_A$ is a weak equivalence for $A$ any positive cofibrant 
symmetric spectrum.  Given a cofibrant symmetric spectrum $B$, choose a 
positive cofibrant replacement $\phi \mc cB \to B$.  
Since $\mI \sm -$ preserves weak equivalences between cofibrant
objects and positive cofibrant objects are cofibrant, $\mI \sm cB \to
\mI \sm B$ is a weak equivalence.  Then one can choose fibrant replacements
and a lift $(\mI \sm \phi)^f$ so that $\psi_B \circ \phi = 
\Hom(\mI, (\mI \sm \phi)^f) \circ
\psi_{cB}$.   Thus, $\psi_B$ is a weak equivalence,
since $\Hom(\mI, -)$ preserves weak equivalences between
fibrant objects and $\phi$ and $\psi_{cB}$ are weak
equivalences. 

If $\mI$ is cofibrant then $\Hom(\mI, -)$ is a right Quillen adjoint functor
to the usual stable model structures.   So the last statement follows 
similarly.
\end{proof}

For the commutative algebra case we need several more assumptions.
These assumptions are satisfied in each of the known examples of equivalences
of commutative algebra spectra~\cite[\S 16]{mmss} 
and~\cite[Theorem 5.1]{sch-comparison}. 
Let $PX= \vee_{i\geq 0} X^{(i)}/ \Sigma_i$ be the monad on $\C$ 
which defines the commutative ring objects (or more properly, {\em monoids}) 
in $\C$.  Here $X^{(i)}$
denotes the $i$th smash power with $\Sigma_i$ permuting the factors
and $X^{(0)}= \mI$.   
To ensure that
$\Hom(\mI,-)$ is a right Quillen adjoint here we require the weak equivalences
and fibrations in the model category of commutative monoids in $\C$  
to be maps which are underlying weak equivalences or fibrations
in $\C$. 
Unlike associative algebras and modules, criteria for the existence
of such a model category on commutative monoids in $\C$ do not exist 
in the current literature.  
Another one of our assumptions here is that the quotient map
from the extended power to the symmetric power, 
$\Phi_i\mc E\Sigma_i \sm_{Sigma_i} X^{(i)} \to X^{(i)}/\Sigma_i$,  
is a weak equivalence for cofibrant objects $X$ in $\C$.  Since
the monad $P$ does not necessarily preserve weak equivalences, this is 
likely to be one of the criteria required for constructing a 
model category on commutative monoids.  

\begin{hypotheses}\label{hyp-ca}
Let $\C$ be a stable simplicial monoidal
model category such that
\begin{enumerate}
\item for any commutative ring $R'$ in $\C$ the commutative $R'$-algebras in 
$\C$ form a model category with a fibrant replacement functor and with 
fibrations and weak equivalences the underlying fibrations and weak 
equivalences in $\C$,
\item 
$\Phi_i \mc E\Sigma_i \sm_{Sigma_i} X^{(i)} \to X^{(i)}/\Sigma_i$  
is a weak equivalence for cofibrant objects $X$ in $\C$ and 
\item there is a monoidal Quillen equivalence from $\spec$ to $\C$
(or any other equivalent condition from Theorem~\ref{cor-simp}).
\end{enumerate}
\end{hypotheses}

These hypotheses hold for the positive stable model categories on orthogonal 
spectra and symmetric spectra by~\cite[10.4, 15.1, 15.2 and 15.5]{mmss} and 
they hold for $S$-modules by~\cite[III.5.1, VI.4.8]{ekmm} 
and~\cite{sch-comparison}.  
These hypotheses are more subtle than those for modules and algebras; for 
example, the first hypothesis does not hold for the usual stable 
model category on symmetric spectra; see~\cite[Section 14]{mmss}.  

\begin{theorem}\label{thm-ca}
Assume $\C$ satisfies Hypotheses~\ref{hyp-ca}. 
Then $\mI \sm -$ and $\Hom(\mI, -)$ induce a Quillen equivalence 
\begin{enumerate}
\item 
from the commutative symmetric ring spectra to 
the commutative rings in $\C$ and 
\item 
from the commutative $R$-algebras for $R$ a cofibrant commutative symmetric
ring spectrum to 
commutative $(\mI \sm R)$-algebras. 
\end{enumerate}
\end{theorem}

\begin{proof}
The first statement is a special case of the
second with $R=\mS$ and $\mI \sm \mS \iso \mI$.  
Since the weak equivalences and fibrations are determined on the
underlying category the restriction of $\Hom(\mI, -)$ is a
right Quillen adjoint functor.  Since the equivalent conditions
in Theorem~\ref{cor-simp} hold, $\Hom(\mI, -)$ preserves
and detects weak equivalences between fibrant objects.  
By~\cite[Lemma 4.1.7]{hss} it thus suffices to show that $\psi_A\mc
A \to \Hom(\mI, (\mI \sm A)^f)$ is a weak equivalence for cofibrant commutative
$R$-algebras where $(-)^f$ is the given fibrant replacement functor.  
Since fibrations are determined on
the underlying category a fibrant replacement as a commutative $R$-algebra
restricts to a fibrant replacement in $\C$.  
We first consider the case where $R=\mS$, the sphere spectrum. 

Let $A= PX$ for a positive cofibrant symmetric spectrum $X$.
We claim that to show $\psi_A$ is a weak equivalence it suffices to show that 
$\psi$ is a weak equivalence for each 
symmetric power $X^{(i)}/\Sigma_i$.  $\Hom(\mI, (\mI \sm -)^f)$ does not 
necessarily
commute with coproducts, but it does commute up to weak equivalence
with homotopy coproducts because it 
is naturally isomorphic to the identity on $\Ho(\spec)$. 
Here the coproduct is a homotopy
coproduct because it is created levelwise and each level of each summand is 
cofibrant.  This is our general strategy; we follow the outline
of the proof of~\cite[Theorem 0.7]{mmss} but there the composite of the 
adjoints commutes 
with colimits on the nose and here it may only commute up to weak
equivalence with homotopy colimits. 
But each of the colimits we must consider is in fact a homotopy colimit of 
symmetric spectra because such homotopy colimits can be computed 
levelwise~\cite[2.2.1]{thh}. 

Now consider each summand.  Since $\mI \sm -$ is a left Quillen adjoint on
the spectrum level it commutes with colimits and smash products with
spaces, so 
$\mI \sm -$ applied to $\Phi_i$ of $X$ in $\spec$ is isomorphic to the
weak equivalence $\Phi_i$ of the cofibrant object $\mI \sm X$ in $\C$.  
Thus it is sufficient to show that $\psi$ is a weak 
equivalence for $E\Sigma_{i+} \sm_{\Sigma_i} X^{(i)}$  
because $\Phi_i$ in $\spec$ and $\Hom(\mI, (\mI \sm \Phi_i)^f)$ are weak 
equivalences by~\cite[15.5]{mmss} and the second hypothesis on $\C$.  
The extended power is the homotopy colimit of the $\Sigma_i$ action on 
$X^{(i)}$.
So since $\psi$ is a weak equivalence for the positive cofibrant
object $X^{(i)}$ and $\Hom(\mI, (\mI \sm -)^f)$ commutes with
homotopy colimits, we conclude that $\psi$ is a weak equivalence
on the extended power.  

Following~\cite[15.9]{mmss}, we proceed by building cofibrant $\mS$-algebras  
using modified generating cofibrations.  Let $\Delta[n]$ denote the
simplicial $n$-simplex and $\dot{\Delta}[n]$ its simplicial boundary.
Let $B_{*}(K_+, K_+, K_+)$ be the simplicial bar construction which is the 
bisimplicial set with $s,t$-simplices $(K_s \times (\Delta[1])_t)_+$.  
Its geometric realization is isomorphic to $K_+ \sm \Delta[1]_+$. 
The inclusion $i_1\mc \Delta[0]_+ \to \Delta[1]_+$ induces
an inclusion of the horizontally constant bisimplicial set 
$c_*(K_+)$ into $B_*(K_+, K_+, K_+)$.  Define $B_*(K_+, K_+, S^0)$ as the pushout
of this inclusion over the map $K_+ \to S^0$.  The geometric realization
of the composite gives a map $K_+ \to B(K_+,K_+,S^0)$ with the
geometric realization $B(K_+,K_+,S^0)$ isomorphic to the unreduced cone
$(CK)_+$.   We can use these composite maps with $K=\dot{\Delta}[r]$ 
instead of the 
simplicially homotopic maps $\dot{\Delta}[r]_+ \to \Delta[r]_+$
to construct generating cofibrations.      
The model category of commutative symmetric ring spectra is cofibrantly
generated with $PF^+I = \{P(F_n\dot{\Delta}[k]_+) \to P(F_n (C\dot{\Delta}[k])_+)
\}_{k\geq 0, n > 0}$ a set of generating cofibrations.     

So we need to show that $\psi_A$ is a weak equivalence when $A$ is a 
$PF^+I$-cell complex~\cite[2.1.18]{hovey-book}.  We have shown $\psi_A$ is a weak equivalence when $A$ is built in one stage.
We next consider $\psi_A$ for $A$ constructed from $PF^+I$ by finitely many pushouts.
Assume the result for those $\mS$-algebras built in $n$ stages, and consider
$A= A_n \sm_{PX} PY$ with $A_n$ built in $n$ stages and $X \to Y$ a coproduct of
maps in $PF^+I$.  Since $F_n$ commutes with colimits and smash products with spaces, it
commutes with geometric realization and the bar construction above.  If $X=\vee_i F_{n_i}
\dot{\Delta}[r_i]_+$ and $T= \vee_i F_{n_i}S^0$, then 
$Y\iso B(X, X, T)$, the geometric realization of the simplicial symmetric spectrum
with $q$-simplices the coproduct of $q+1$ copies of $X$ and one copy of $T$.   
Statements analogous to~\cite[5.1]{mmss} and~\cite[VII.2.10, VII.3.3]{ekmm} show that
the category of commutative $\mS$-algebras is tensored over simplicial
sets and the underlying symmetric spectrum of the geometric realization of a 
simplicial commutative $\mS$-algebra is isomorphic to the geometric realization 
of the underlying simplicial symmetric spectrum.
Since $P$ commutes with colimits and converts smash products with spaces to
tensors with spaces, $P$ commutes with geometric realizations.  Hence,     
\[A\iso A_n \sm_{PX} PY \iso A_n \sm_{PX} B(PX, PX, PT) \iso B(A_n, PX, PT).
\]     

Since tensors with simplicial sets and colimits in $\spec$ are levelwise,
the geometric realization is constructed on each level.   
But the geometric realization of a bisimplicial
set is weakly equivalent to the homotopy colimit by~\cite[XII.4.3]{BK}.   
So the 
geometric realization $B(A_n, PX, PT)$ is weakly equivalent to the homotopy 
colimit and $\Hom(\mI, (\mI \sm -)^f)$ commutes with the homotopy colimit. 
So it is enough
to show that $\psi$ is a weak equivalence on each simplicial level of 
$B(A_n, PX, PT)$.
The $q$-simplices here are given by $A_n \sm (PX)^{(q)} \sm PT \iso 
A_n \sm P(X \vee \cdots \vee X \vee T)$.  These $q$-simplices can be constructed
in $n$ stages, so by induction $\psi$ is a weak equivalence here, as required.

Finally, we must consider a filtered colimit of these commutative $\mS$-algebras built
in finitely many stages.  
Filtered colimits of commutative $\mS$-algebras are created on the
underlying symmetric spectra which in turn are created on each level.  Since the
maps in question here are constructed from $PF^+I$, they are injections.  
So the colimit is over level injections between level cofibrant objects and it 
is weakly equivalent to the homotopy colimit.  By induction we have shown that
$\psi$ is a weak equivalence at each spot in the colimit and 
$\Hom(\mI, (\mI \sm -)^f)$
commutes with homotopy colimits, so $\psi$ is a weak equivalence on the colimit as well.  
So we conclude that $\psi_A$ is a weak equivalence for any cofibrant 
commutative symmetric ring spectrum $A$ (i.e. any retract of a
$PF^+I$-cell complex). 

For the second statement we consider cofibrant commutative $R$-algebras $A$ 
for $R$
a cofibrant commutative symmetric ring spectrum.  Since $R$ is cofibrant the
unit map $\mS \to R$ is a cofibration of commutative symmetric ring spectra. Since $A$ is 
cofibrant as a commutative $R$-algebra, the unit map $R \to A$ and hence the
composite $\mS \to A$ is also a cofibration.  Hence, by the above, $\psi_A$ is a
weak equivalence.  
\end{proof}

\section{Non-simplicial case}\label{sec-non-simp}

In this section we consider the case when the given stable, monoidal
model category $\C$ is not simplicial.  
We produce a Quillen adjoint pair from $\spec$
to $\C$ whose derived functors are monoidal.   This can be used for example
to produce a Quillen adjoint pair from $\spec$ to $\mZ$-graded chain complexes. 
Since $\C$ is not simplicial, we need a new definition for 
a desuspension of the unit.  

\begin{definition}\label{def-non-simp}
{\em
A {\em cylinder object} for a cofibrant object $X$
is an object $X \times I$ with a factorization of the fold map
$X \amalg X \varr{i} X \times I \varr{p} X$ such that $i$ is a cofibration
and $p$ is a weak equivalence.  A model for the {\em suspension}, $\Sigma X$,
is the cofiber of $X\amalg X \varr{i} X \times I$ for some cylinder 
$X \times I$~\cite[Chapter I \S 1, 2]{Q}.  A {\em good desuspension} of the unit is a cofibrant object
$\mIo$ with a weak equivalence $\eta \mc \Sigma \mIo \to \mI$ for some model
of the suspension.  A {\em stable monoidal model category} is a monoidal
model category which is stable and has a good desuspension of the unit.
} 
\end{definition}

\begin{theorem}\label{thm-non-simp}
Let $\C$ be a stable monoidal model category. 
Then there is a Quillen adjoint pair from  
the positive stable model structure on $\spec$ to $\C$, again denoted 
by $\mI \sm -$ and $\Hom(\mI, -)$, such that the total left derived
functor $\mI \sm^L -$ is strong monoidal.  Moreover, $\Hom(\mI, -)$ is lax
monoidal, $\mI \sm \mS \iso \mI$, and $\mI \sm Q\mS \to \mI \sm \mS$
is a weak equivalence.
\end{theorem}

As with the cofibrant desuspension, the existence of a functor with
the properties listed for $\mI \sm -$ implies the existence of a good
desuspension.  
Under additional conditions on the unit
this monoidal functor induces a monoidal equivalence of
the homotopy category of $\C$ and the homotopy category of symmetric
spectra.  The proof of the following statement is similar 
to Theorem~\ref{cor-simp}. 
\begin{theorem}\label{cor-non-simp}
Let $\C$ be a stable monoidal  model category. 
The following conditions are equivalent:
\begin{enumerate}
\item There is a $\pi^s_*$-linear, triangulated equivalence between
the homotopy category of $\C$ and the homotopy category of $\spec$ which
takes the unit $\mI$ of the monoidal product to the unit 
$\mS$. 
\item $\mI$, is a small weak generator and  
$[\mI, \mI]_*^{\Ho(\C)}$ is freely generated as a $\pi_*^s$-module by 
$id_{\mI}$.
\item $\C$ and $\spec$ are Quillen equivalent via functors whose derived
functors are strong monoidal.
\item There is a $\pi^s_*$-linear, triangulated, monoidal equivalence 
between $\Ho(\C)$ and $\Ho(\spec)$. 
\end{enumerate}
\end{theorem}

To construct the right adjoint $\Hom(\mI, -)$ we use {\em cosimplicial
resolutions} since $\C$ is not simplicial.   These 
were first used in~\cite{DK} to construct function complexes on 
homotopy categories, but in~\cite{DHK} this theory has been extended to provide
function complexes on model categories.   Our main reference here 
is~\cite[Chapter 5]{hovey-book}, see also~\cite{SS3}.  

Given a cosimplicial object $X^{\bd}$ in $\C^{\Delta}$ and a pointed 
simplicial set $K$ denote 
the coend~\cite[Chapter IX \S 6]{maclane} in $\C$ by $X^{\bd} \otimes_{\Delta} K$ see 
also~\cite[Chapter 5 \S 7]{hovey-book}.  Define $X^{\bd}\otimes K$ by $(X^{\bd} \otimes K)^n
= X^{\bd} \otimes_{\Delta} (K \sm \Delta[n]_+).$
Notice $X^{\bd} \otimes K$ and $X^{\bd} \otimes_{\Delta} K$ are objects in 
different categories ($\C^{\Delta}$ and $\C$ respectively.)  Set $S^m =
(S^1)^m$ and denote $X^{\bd} \otimes S^m$ by $\Sigma^m(X^{\bd})$.
If $X^{\bd}$ is a 
cosimplicial object and $Y$ is an object of $\C$ then $\C(X^{\bd}, Y)$ 
is a simplicial set with degree $n$ the set of $\C$-morphisms $\C(X^n, Y)$.
There is an adjunction isomorphism $\C(X^{\bd} \otimes K, Y) \iso 
\map(K,\C( X^{\bd}, Y))$.  
This shows that $X^{\bd} \otimes (K \sm L)
\iso (X^{\bd} \otimes K) \otimes L$ since they both represent the same functor.
In particular, $\Sigma^m(X^{\bd})$ is the $m$th iterated suspension of 
$X^{\bd}$.

We consider the Reedy model category on
$\C^{\Delta}$, the cosimplicial objects on $\C$~\cite{Reedy},
\cite[Theorem 5.2.5]{hovey-book}.  An object $A^{\bd}$ is {\em Reedy cofibrant}
if the map $A^{\bd} \otimes \partial\Delta[k]_+ \to  A^{\bd} \otimes 
\Delta[k]_+ \iso A^k$ 
is a cofibration for each $k$.  A {\em cosimplicial resolution}
is then a Reedy cofibrant object of $\C^{\Delta}$ such that each of
the codegeneracy and coface maps are weak equivalences.  That is, cosimplicial
resolutions are the
Reedy cofibrant, homotopically constant objects.   A cosimplicial resolution
$A$ is called a cosimplicial frame of the cofibrant object $A^0$ 
in~\cite[Chapter 5]{hovey-book}.  

The category $\C^{\Delta}$ has a symmetric monoidal product, defined
on each level by the symmetric monoidal product on $\C$.  That is, 
$(A^{\bd} \sm B^{\bd})^n \iso A^n \sm B^n$.   
The following proposition collects several useful properties of the 
proceeding constructions.  These properties follow 
from~\cite[Theorem 16.4.2, Proposition 16.11.1]{psh}
and~\cite[Proposition 5.7.1 and 5.7.2]{hovey-book}. 

\begin{proposition}\label{prop-frames}
Let $A^{\bd}$ and $B^{\bd}$ be cosimplicial resolutions in $\C^{\Delta}$
where $\C$ is a monoidal model category such that $\sm$ commutes with colimits.
\begin{enumerate}
\item \label{one} $A^{\bd} \sm B^{\bd}$ is a cosimplicial resolution.
\item \label{two} $\Sigma A^{\bd}$ is a cosimplicial resolution.
\item \label{three} There is a natural level
equivalence $\Sigma (A^{\bd} \sm B^{\bd}) \to (\Sigma A^{\bd}) \sm B^{\bd}.$
\end{enumerate}
\end{proposition}

\begin{proof}
For part~\ref{one}, note that $A^{\bd} \sm B^{\bd}$ is homotopically constant 
since the smash product of two maps which are each weak equivalences between 
cofibrant objects is a weak equivalence.
The monoidal product also preserves Reedy cofibrant objects. 
By~\cite[Proposition 16.11.1]{psh}, since $\sm \mc \C \times \C \to \C$ preserves 
cofibrations  and $\sm$ commutes with colimits, 
the prolongation $\sm \mc \C^{\Delta} \times \C^{\Delta} \to \C^{\Delta}$
also preserves cofibrations.  This uses~\cite[Theorem 16.4.2]{psh} to recognize
that the Reedy model category on $\C^{\Delta} \times \C^{\Delta}$ agrees with
the Reedy model category on $(\C \times \C)^{\Delta}$.

For part~\ref{two}, the map $\Sigma A^{\bd} \otimes 
(\partial\Delta[k]_+ \to \Delta[k]_+)$ is 
isomorphic to the map $A^{\bd} \otimes (S^1 \sm \partial\Delta[k]_+ \to
S^1 \sm \Delta[k]_+)$.   By~\cite[Proposition 5.7.1]{hovey-book}, if $A^{\bd}$ 
is Reedy
cofibrant then this map is a cofibration.  So $\Sigma A^{\bd}$ is Reedy
cofibrant.  Since each map $S^1 \sm \Delta[n]_+ \to S^1 \sm \Delta[n+1]_+$
is a trivial cofibration, the coface maps of $\Sigma A^{\bd}$ are trivial 
cofibrations~\cite[Proposition 5.7.2]{hovey-book}.  Since $s^i d^i = \id$
the codegeneracy maps are also weak equivalences.

For part~\ref{three}, the coend defining $\Sigma(A^{\bd} \sm B^{\bd})^m$ is a colimit of
copies of $(A^k \sm B^k)$ indexed by the non-base point $k$-simplices of 
$S^1 \sm \Delta[m]_+$.  Use the map $k \to m$ in $\Delta$ determined by
the $k$-simplex of $\Delta[m]_+$ to induce  
a map $B^k \to B^m$.  These maps are all compatible and define a map 
$\Sigma(A^{\bd} \sm B^{\bd})^m \to \Sigma(A^{\bd} \sm B^m)^m \iso 
\Sigma(A^{\bd})^m \sm B^m \iso (\Sigma(A^{\bd}) \sm B^{\bd})^m$.  Here
$(A^{\bd} \sm B^m)^k \iso A^k \sm B^m$ and 
we have used the fact that $\sm$ commutes with colimits.
Since parts~\ref{one} and~\ref{two} 
show that this is a map between cosimplicial frames, it is a level equivalence
if degree zero is a weak equivalence.  The map $A^{\bd} \sm B^{\bd} \to
A^{\bd} \sm B^0$ is a level equivalence between Reedy cofibrant objects
by part~\ref{one} and the monoidal model structure on $\C$.
Hence $(A^{\bd} \sm B^{\bd}) \otimes_{\Delta} S^1 \to
(A^{\bd} \sm B^0)\otimes_{\Delta} S^1 \iso (A^{\bd} \otimes_{\Delta} S^1) \sm B^0$ is a weak 
equivalence by~\cite[Proposition 5.7.1]{hovey-book}.
\end{proof}

\begin{proof}[Proof of Theorem~\ref{thm-non-simp}]
To define the right adjoint $\Hom(\mI, -)$ we consider cosimplicial 
resolutions 
related to $\mI$.  First, let $\omega^0 \mI$ be the constant cosimplicial 
object on $\mI$.  Since $\mI$ is not necessarily cofibrant, $\ozI$
is not necessarily a cosimplicial resolution.
Since $\C$ has a good desuspension of the unit,  one can build a cosimplicial 
resolution $\ooI$ of the cofibrant object $\mIo$.  
Define $(\ooI)^0=\mIo$ and $(\ooI)^1= \mIo \times I$.  Define the coface maps 
as the two inclusions $X \to X \amalg X \varr{i} X\times I$ and define the
codegeneracy map as the map $X \times I \varr{p} X$.  Using the factorization
properties in $\C$ one can inductively define the higher levels of $\ooI$,
see the proof of~\cite[Theorem 5.1.3]{hovey-book}.   Since $\Sigma(X^{\bd})$
is the cofiber of $X^{\bd} \otimes (S^0 \vee S^0) \to X^{\bd} \otimes 
\Delta[1]_+$, $\Sigma(\ooI)^0$ is the cofiber of $i$, that is, 
a model for the suspension of $\mIo$.  
Then the weak equivalence $\eta\mc \Sigma\mIo \to \mI$
extends to a level equivalence $\eta^{\bd} \mc \Sigma \ooI \to \ozI$.  
Define $\onI = (\ooI)^{\sm n}$ for $n > 0$.  
Define the right adjoint $\Hom(\mI, -)$ in level $n$ to
be $\C(\onI, -)$.  The symmetric group on $n$ letters 
permutes the factors of $\onI$.   The structure maps are induced by 
the map $\eta$.   
Proposition~\ref{prop-frames}
part~\ref{three} provides a level equivalence $\phi\mc \Sigma^m(A^{\sm m} \sm B)
\to (\Sigma A)^{\sm m} \sm B$.
The isomorphism of $\Sigma^m(X^{\bd})$ with
the $m$-fold iterated suspension of $X^{\bd}$ induces a 
$\Sigma_m \times \Sigma_n$-equivariant level equivalence where $\Sigma_m$
acts trivially on the target:
\[ \Sigma^m(\ompnI) \varr{\phi} (\Sigma \ooI)^{\sm m} \sm \onI
\varr{(\eta^{\sm m}) \sm \id}(\ozI)^{\sm m} \sm \onI \iso \onI. \]  
Applying $\C(-, Z)$ to the displayed composition and taking adjoints gives the 
$\Sigma_m \times \Sigma_n$-equivariant structure map
\[ S^m \sm \C(\onI, Z) \to \C(\ompnI, Z). \] 

Let $Q\ozI$ denote the constant cosimplicial object on $Q\mI$,
the chosen cofibrant replacement of $\mI$.  Then since degree zero of $\eta$
factors through $Q\mI$, $\eta$ factors as two level equivalences
$\Sigma \ooI \to Q\ozI \to \ozI$.
Hence $(\eta)^{\sm m} \sm \id_{A^{\bd}}$ for any cosimplicial resolution 
$A^{\bd}$ 
is a level equivalence by the monoidal model structure on $\C$.  
Since $\onI$ for $n>0$ is a cosimplicial resolution by Proposition~\ref{prop-frames}
part~\ref{one}, each map $\Sigma^m(\ompnI) \to \onI)$ with $n> 0$ is a weak 
equivalence.  By the pointed version of~\cite[Corollary 5.4.4]{hovey-book}, 
$\C(A^{\bd}, -)$ preserves fibrations and trivial fibrations when
$A^{\bd}$ is a cosimplicial resolution and for $Z$ fibrant $\C(-, Z)$ takes
level equivalences between cosimplicial resolution to weak equivalences.    
So $\Hom(\mI, -)$ takes fibrant objects to positive $\Omega$-spectra 
and (trivial)
fibrations to positive level (trivial) fibrations.  Thus $\Hom(\mI, -)$
is a right Quillen functor by~\cite[Corollary A.2]{dugger} since positive 
stable fibrations between positive $\Omega$-spectra are positive level 
fibrations by Proposition~\ref{prop-pos-stable}.  The left adjoint,
$\mI \sm -$ is formed as in the simplicial case except here the tensors
of cosimplicial resolutions with simplicial sets are given by coends. 

To show that the total left derived functor $\mI \sm^L -$ is strong
monoidal we first show that $\Hom(\mI, -)$ is 
lax monoidal.  For the unit map take the non-base point of $S^0$ to the
identity map in simplicial degree zero of $\Hom(\mI, \mI)^0$.  
The monoidal product on $\C$ induces a natural
map $\C(\omI, A) \sm \C(\onI, B) 
\to \C(\ompnI, A \sm B)$. Assembling these levels produces a 
natural map $\Hom(\mI, A) \sm \Hom(\mI, B) \to \Hom(\mI, A \sm B)$. 
Hence, $\Hom(\mI, -)$ is lax monoidal. 
So its left adjoint
$\mI \sm -$ is lax comonoidal.  Also, $\mI \sm F_n'X \iso (\onI 
\otimes_{\Sigma_n} X)^0$ because they represent the same functor in $\C$.  
So $\mI \sm -$ takes 
the cofibrant replacement $Q\mS \iso F_1'S^1 \to F_0S^0 \iso \mS$ 
to the weak equivalence
$\eta \mc (\ooI \otimes S^1)^0 \to \mI$. 
The comonoidal structure on $\mI \sm -$ induces a natural transformation
$\mI \sm^L (A \sm^L B) \to (\mI \sm^L A) \sm^L (\mI \sm^L B)$
where $\mI \sm^L -$ is the total left derived functor of $\mI \sm -$.   
Since $\mI \sm^L \mS \iso \mI$, 
this map is an isomorphism for $A = \mS$ and any
$B$.  For fixed $B$ both the source and target are exact functors in 
$A$ which commute with coproducts,  so the  
objects $A$ where this transformation is an isomorphism form a localizing
subcategory which contains the generator $\mS$.  Hence this transformation
is an isomorphism for all $A$ and $B$.  So $\mI \sm^L -$ is strong
monoidal.
\end{proof}

\begin{remark}\label{rem-non-simp-desuspension}
{\em Let $\C$ be a monoidal model category with a Quillen adjoint pair between 
$\C$ and the positive stable model category on $\spec$ with left adjoint
$L\mc \spec \to \C$.  If $L(Q\mS) \to L(\mS)$ is a weak equivalence and 
$L(\mS)\iso \mI$ then $\C$ has
a good desuspension of the unit.  Define $\mIo = L(F_1 S^0)$, with
cylinder $L(F_1 \Delta[1]_+)$ and model of the suspension $L(F_1 S^1)$.
These definitions have all the necessary properties since 
$L$ preserves positive cofibrations and weak equivalences between positive
cofibrant objects.  
Define $\eta \mc L(F_1S^1) \to L(F_0 S^0)$ as the adjoint of the identity 
map on level one.  
In fact, the cosimplicial resolution $\ooI$ can be defined 
by $(\ooI)^n = L(F_1\Delta[n]_+)$.  
}\end{remark}

\end{document}